\documentclass[12pt]{amsart}

\usepackage{amsthm}
\usepackage{enumerate}
\usepackage{amsmath}
\usepackage{mathrsfs}

\usepackage[colorlinks]{hyperref}

\newtheorem{thm}{Theorem}[section]
\newtheorem{lem}[thm]{Lemma}
\newtheorem{prop}[thm]{Proposition}
\newtheorem{cor}[thm]{Corollary}
\newtheorem*{nota}{Notation}

\newtheorem{defn}[thm]{Definition}

\title[Invariant Means on $VN^n(G)$]{Invariant Means on $VN^n(G)$}
\author{Kanupriya}
\address{Kanupriya,\newline\indent Department of Mathematics,\newline\indent Indian Institute of Technology Delhi,\newline\indent New Delhi - 110016, India.}
\email{kanupriyawadhawan3@gmail.com}
\author{N. Shravan Kumar}
\address{N. Shravan Kumar,\newline\indent Department of Mathematics,\newline\indent Indian Institute of Technology Delhi,\newline\indent New Delhi - 110016, India.}
\email{shravankumar.nageswaran@gmail.com}

\begin{document}

\begin{abstract}
    Let $G$ be a locally compact group, and $VN^n(G)$ is the dual of the multidimensional Fourier algebra $A^n(G)$. In this article, we define invariant means on $VN^n(G)$ and prove that the set of all invariant means on $VN^n(G)$ is non-empty. Further, we investigated the invariant means on $VN^n(G)$ for discrete and non-discrete cases of $G$. Also, we show that if $H$ is an open subgroup of $G$, then the number of invariant means on $VN^n(H)$ is the same as that of $VN^n(G)$. Finally, we study invariant means on the dual of the algebra $A_0^n(G)$, the closure of Fourier algebra $A^n(G)$ in the cb-multiplier norm.
\end{abstract}
\maketitle
\section{Introduction}
Invariant means on $VN(G)$ are important in studying group von Neumann algebras and applications in functional analysis, operator algebras, and geometric group theory. An invariant mean on $VN(G)$ is a state (a positive linear functional of norm $1$) that remains invariant under the action of $G$. In 1972, Renaud \cite{R} developed and analysed invariant means on $VN(G)$. Further, in 1982, C. Chou \cite{C} investigated the links between invariant means and other algebraic features of von Neumann algebras. One of the main results of  Renaud's paper is that a unique invariant mean on $VN(G)$ exists if and only if $G$ is discrete. Later, Lau and Losert \cite{LL} continue to study the invariant means for the case when $G$ is non-discrete. Since then, many authors have attempted to generalise Renaud's classical insights. For the case of homogeneous spaces, Cho-Ho Chu and A. T. M. Lau \cite{CL} have extended this. N. Shravan Kumar \cite{S}, initiated a study of invariant means on $VN(K)$, where $K$ is the ultraspherical hypergroup of $G$. Further, the author studied some generalised translations and generalised invariant means in \cite{S2}. This article focuses on the study of the invariant means on $VN^n(G)$.
 
 The multidimensional Fourier algebra $A^n(G)$ on a locally compact group $G$ was introduced by Todorov and Turowska \cite{ToTu}, where the authors studied multipliers on the multidimensional Fourier algebra. The idea of a multidimensional version was motivated by the fact that, when $G$ is abelian, the two-dimensional Fourier-Stieltjes algebra, via the Fourier transform, coincides with the space of bimeasures on $\widehat{G}\times\widehat{G}$ \cite{GrSc}. The dual space of the multidimensional Fourier algebra $A^n(G)$ is $VN^n(G)$ which is defined in the same paper \cite{ToTu}. We investigate the invariant means on this dual algebra $VN^n(G)$. 
 
 In Section \ref{Sctn3CS5}, we first introduced the notion of means on the dual algebra $VN^n(G)$ and studied some properties of the spaces $M_{A^n(G)}$ and $M_{B^n(G)}$. In the next section \ref{sctn4CS5}, the concept of invariant means is discussed and
subsequently noticed that identical to the case of $VN(G)$, the set of all invariant means on $VN^n(G)$, denoted by $\operatorname{TIM}^n(\widehat{G})$, is a non-empty set. One of the key findings here is that the set $\operatorname{TIM}^n(\widehat{G})$ is a singleton set when the group $G$ is discrete. 
 An important result of this study also appears in this section \ref{sctn4CS5}, which deals with the case when $G$ is non-discrete. The existence of more than one invariant means on $VN^n(G)$ is proved for the case when $G$ is non-discrete. Further, in section \ref{sctn5CS5}, we establish a relationship between the invariant means on $VN^n(H)$ and the invariant means on $VN^n(G)$ for an open subgroup $H$ of $G$. Particularly, we prove that the total number of invariant means on $VN^n(H)$ is equal to the total number of invariant means on $VN^n(G)$.

Forrest and Miao studied invariant means in the algebra $A_0(G)$ \cite{FM}. Here, $A_0(G)$ represents the closure of $A(G)$ in the cb-multiplier norm. In the last section \ref{sctn6CS5}, we examined invariant means on the dual of the algebra $A^n_0(G)$, which is the closure of the multidimensional Fourier algebra $A^n(G)$ in the cb-multiplier norm.

We shall now begin with some preliminaries that are needed in the sequel.
\section{Preliminaries and Notations}
        This section provides a comprehensive overview of fundamental concepts in operator spaces. Given the significant relevance of tensor products, the primary objective of this part is to get the necessary background information on these subjects.

        Let $X$ be a vector space. The set of all $n\times n$ matrices with elements in the space $X$ will be represented as $M_n(X)$. An operator space is a complex vector space $X$ equipped with a norm $\|\cdot\|_n$ on the matrix space $M_n(X)$ for every natural number $n$, satisfying the following conditions:
        \begin{enumerate}[(i)]
            \item $\|x\oplus y\|_{m+n}=max\{\|x\|_m,\|y\|_n\}$ and
            \item $\|\alpha x\beta\|_n\leq\|\alpha\|\|x\|_m\|\beta\|$
        \end{enumerate}
        for all $x\in M_m(X),$ $y\in M_n(X),$ $\alpha\in M_{n,m}$ and $\beta\in M_{m,n}.$ 

        Let $X$ and $Y$ be operator spaces, and let $\varphi:X\rightarrow Y$ be a linear map. The $n^{\text{th}}$-amplification of the linear transformation $\varphi,$ abbreviated $\varphi_n,$ is defined as a mapping from the matrix space $M_n(X)$ to $M_n(Y),$ where $n$ is a natural number. This mapping is given by $\varphi_n([x_{ij}])=[\varphi(x_{ij})]$. The linear transformation $\varphi$ is said to be completely bounded if $\sup\{\|\varphi_n\|:n\in\mathbb{N}\}<\infty.$ We shall use the notation $\mathcal{CB}(X,Y)$ to represent the set of all completely bounded linear mappings from $X$ to $Y$, with the following norm: $\|\cdot\|_{cb},$ $$\|\varphi\|_{cb}:=\sup\{\|\varphi_n\|:n\in\mathbb{N}\},\ \varphi\in\mathcal{CB}(X,Y).$$ We shall say that $\varphi$ is a {\it complete isometry (complete contraction)} if $\varphi_n$ is an isometry (a contraction) for each $n\in\mathbb{N}.$

        Ruan's theorem implies that for each abstract operator space $X$, there exists a Hilbert space $\mathcal{H}$ and a closed subspace $Y$ of $\mathcal{B}(\mathcal{H})$ such that $X$ and $Y$ are completely isometric. 

        For operator spaces $X$ and $Y$ the {\it Haagerup tensor norm} of $u\in M_n(X\otimes Y)$ is given as $$\|u\|_h=\inf\left\{ \|x\|\|y\|:u=x\odot y,x\in M_{n,r}(X), y\in M_{r,n}(Y),r\in\mathbb{N} \right\}.$$ 
        The expression $x \odot y$ represents the inner matrix product of $x$ and $y$. This product is defined as $(x \odot y)_{i,j} = \underset{k=1}{\overset{r}{\sum}} x_{i,k} \otimes y_{k,j}$, where $x \odot y$ belongs to the matrix space $M_{n}(X \otimes Y)$. Refer to \cite[Chapter 9]{ER3} for more details. The norm $\|\cdot\|_h$ is an operator space norm.  We shall use the notation $X\otimes^h Y$ for the resulting operator space which is obtained by completion of $X \otimes Y$ with respect to the Haagerup norm. If $X$ and $Y$ are C*-algebras and $n=1$, the Haagerup norm may be expressed as follows. For $u\in X\otimes^h Y,$ 
        $$\|u\|_h=\inf \left\{ \left\| \underset{n\in\mathbb{N}}{\sum} x_n x_n^\ast \right\|_X^{1/2} \left\| \underset{n\in\mathbb{N}}{\sum} y_n^\ast y_n \right\|_Y^{1/2} :u= \underset{n\in\mathbb{N}} {\sum}x_n\otimes y_n \right\}.$$ 
            
        The {\it extended Haagerup tensor product} of $X$ and $Y,$ denoted by $X\otimes^{eh} Y$, is defined as the space of all normal multiplicatively bounded functionals on $X^\ast\times Y^\ast.$ According to \cite{ER2}, $(X\otimes^h Y)^\ast$ and $X^\ast\otimes^{eh}Y^\ast$ are completely isometric. For the dual operator spaces $X^\ast$ and $Y^\ast,$ the {\it $\sigma$-Haagerup tensor product (or normal Haagerup tensor product)} is defined by $$X^\ast\otimes^{\sigma h}Y^\ast=(X\otimes^{eh}Y)^\ast.$$ Also, the following inclusions hold completely isometrically: $$X^\ast\otimes^h Y^\ast \hookrightarrow X^\ast\otimes^{eh} Y^\ast \hookrightarrow X^\ast\otimes^{\sigma h} Y^\ast.$$

       \begin{nota}
            The n-tensor product of a space X is denoted by $$\otimes_{n}X=\underbrace{ X \otimes X \otimes\cdots\otimes X}_n.$$ Now, we define $\otimes_{n}^h X,$ $\otimes_{n}^{eh} X,$ and $\otimes_{n}^{\sigma h} X$ as the completion of $\otimes_{n}X$ with respect to their corresponding norms.
        \end{nota}

       Please refer to \cite{ER1,ER2} to learn more about these tensor products. To obtain more comprehensive information on operator spaces, readers should look into the references \cite{ER3} or \cite{Pis}.

    Let $G$ be a locally compact group and fix a left Haar measure $dx$ on $G$. The symbol $\lambda_G$ represents the left regular representation of the group $G$ acting on the Hilbert space $L^2(G)$ by left translations. The function $\lambda_G(s)f(t)$ is defined as $f(s^{-1}t)$, where $s$ and $t$ belong to the group $G$, and $f$ is in $L^2(G)$. The group von Neumann algebra $VN(G)$ is the smallest self-adjoint subalgebra in $\mathcal{B}(L^2(G))$ that contains all $\lambda_G(s)$ for all $s \in G$ and is closed in the weak operator topology. 
    
    The Fourier algebra $A(G)$ is the predual of the group von Neumann algebra $VN(G)$ and every element $u$ in $A(G)$ is represented by $u(x) = \langle \lambda_G(x)f,g \rangle$ with $\|u\|_{A(G)} = \|f\|_2\|g\|_2.$ The duality between $A(G)$ and $VN(G)$ is given by $\langle T,u \rangle = \langle Tf,g \rangle.$ For further information, readers may refer to the work of Eymard \cite{Eym1}. 
    
    The group algebra $L^1(G)$ is an involutive Banach algebra and its enveloping $C^*$-algebra is denoted by $C^*(G)$. The closure of the left regular representation of $L^1(G)$ in the operator norm is known as the reduced  group $C^*$-algebra, which is denoted by $C_r^*(G)$. We shall denote the group  $G\times G\times\cdots\times G\ (n\mbox{-times})$ by $G^n$.

    A linear functional $m$ on $VN(G)$ is called a mean if it satisfies the following properties:\\
    (i) $m(T) \geq 0\ \forall\ T \geq 0.$\\
    (ii) $m(I) =1$ where $I$ is the  identity operator.\\
     The mean $m$ is called invariant if 
     $$m(uT)= u(e)m(T)\ \forall\ T\in VN(G),\ u \in A(G),$$
      where $uT$ the operator in $VN(G)$ defined by $\langle uT, v\rangle = \langle T, uv\rangle$ for $v \in A(G)$. To know more about the invariant means on $VN(G)$, one can refer to \cite{R}.

    The multidimensional version of Fourier algebra denoted by $A^n(G)$, is actually defined as collection of all functions $f\in L^{\infty}(G^n)$ such that there exists a normal completely bounded multilinear functions $\Phi$ on $\underbrace{VN(G) \times VN(G) \times \cdots \times VN(G)}_n$ satisfying $f(x_1,x_2,...,x_n)= \Phi(\lambda(x_1),...,\lambda(x_n))$. From \cite{ToTu}, $A^n(G)$ coincides with the space $\otimes_n^{eh} A(G).$ Also, the dual of $A^n(G)$ is completely isometrically isomorphic to $VN^n(G):= \otimes_n^{\sigma h}VN(G)$. 
    The multidimensional version of Fourier-Stieltjes algebra is denoted by $B^n(G)$, and is defined as $B^n(G):= \otimes_n^{eh}B(G)$.
    See \cite{ToTu} for more details.
\section{ Means on  $VN^n(G)$ } \label{Sctn3CS5}
\begin{defn}
    A mean $m$ on $VN^n(G)$ is a linear functional on $VN^n(G)$ with the following properties\\
    (i) $m(T) \geq 0\ \forall\ T \geq 0.$\\
    (ii) $m(I\otimes \ldots \otimes I) =1$ where $I$ is identity operator.
\end{defn}
Let $M$ be the set of all means on $VN^n(G)$. One can observe that $M$ is a weak*-compact convex subset of $VN^n(G)^{*}$.

Now, we begin with some notations.
\begin{nota}
 \begin{enumerate}{}
     \item $M_{A^n(G)}:=\left\{u \in A^n(G):\|u\|_{A^n(G)}=u(e,e,\ldots,e)=1\right\}.$ 
     \item $M_{B^n(G)}:=\left\{v \in B^n(G):\|v\|_{B^n(G)}=v(e,e,\ldots,e)=1\right\}.$
     \item For any closed set $E\subseteq G^n$,
     $$j_{A^n(G)}(E) = \{u\in A^n(G):\widehat{u}\mbox{ has compact support disjoint from E} \}.$$
\end{enumerate} 
\end{nota}
\begin{defn}
    \begin{enumerate}
        \item For $m \in( VN^n(G))^*$ and $T \in VN^n(G)$, define $m \circ T \in VN^n(G)$ as
        $$\left\langle m \circ T, v \right\rangle = \left\langle m , vT \right\rangle\ \forall\ v \in A^n(G).$$
        \item For $m,p \in (VN^n(G))^*$, define $m \circ p \in (VN^n(G))^*$ as 
        $$\left \langle m \circ p, T \right\rangle = \left\langle m , p \circ T \right\rangle\ \forall\ T \in VN^n(G).$$
    \end{enumerate}
    \end{defn}
\begin{lem}\label{P1XCS5}
    Let $\widetilde{V}$ be a neighbourhood of $(e,e,\ldots,e)$ in $G^n$. Then there exists a function $u \in A^n(G)$ such that:\\
(a) $0 \leq u \leq 1$;\\
(b) $\|u\|_{A^n(G)}=u(e,e,\ldots,e)=1$;\\
(c) $\operatorname{supp}(u) \subset \widetilde{V}$.
\end{lem}
\begin{proof}
Let $\widetilde{U}$ be a symmetric, relatively compact neighbourhood of $e$ in $G$ such that $\widetilde{U} \times \widetilde{U} \times \ldots \times \widetilde{U} \subset \widetilde{V}$. Then $u=\frac{1}{m_{G}(\widetilde{U})} \chi_{\widetilde{U}} * \chi_{\widetilde{U}} \otimes \frac{1}{m_{G}(\widetilde{U})} \chi_{\widetilde{U}} * \chi_{\widetilde{U}} \otimes \ldots \otimes \frac{1}{m_{G}(\widetilde{U})} \chi_{\widetilde{U}} * \chi_{\widetilde{U}}$ will be the desired function.
\end{proof}

\begin{thm} \label{Thm1XCS5}
    Let $m \in M$ be such that $v \circ m=m$ for all $v \in B^n(G)$, then $v \circ m= v(e,e,\ldots,e) m$.
\end{thm}
\begin{proof}
    \begin{enumerate}[(i)]
        \item Let $v \in B^n(G)$ be such that $v=1$ on a neighbourhood $\widetilde{V}$ of $(e,e,\ldots,e)$ in $G^n$, then consider $u$ to be the function as in Proposition \ref{P1XCS5} for the neighbourhood $\widetilde{V}$. Then $uv=u$. Hence,
        $$v \circ m=v \circ (u \circ m)=(vu) \circ m=u \circ m=m.$$
        \item Let $v \in A^n(G)$ be such that $v(e,e,\ldots,e)=0$. From \cite{Kanu}, $\{(e,e,\ldots,e)\}$ is a set of spectral synthesis and hence there exists a sequence $\left\{v_{n}\right\} \in j_{A^n(G)}(\{(e,e,\ldots,e)\})$ such that $\left\|v_{n}-v\right\|_{A^n(G)} \rightarrow 0$. Also, by using (i), $m=$ $\left(1-v_{n}\right) \circ m=m-v_{n} \circ m$, and hence $v_{n} \circ m=0$. Now,
$$
\|v \circ m\|=\left\|\left(v-v_{n}\right) \circ m\right\| \leq\left\|v-v_{n}\right\|\|m\| \rightarrow 0
$$
Hence, $v\circ m=0.$
\item Let $v \in B^n(G)$ be such that $v(e,e,\ldots,e) \neq 0$. Consider $u \in M_{A^n(G)}$ and let $w \in A^n(G)$ be such that $w=1$ on some neighborhood $\widetilde{V}$ of $(e,e,\ldots,e)$. As $\frac{v u}{v(e,e,\ldots,e)}-w=0$ on $(e,e,\ldots,e)$,\\
by using (ii), $\frac{v u}{v(e,e,\ldots,e)} \circ m=w \circ m$, and by using (i),

\begin{align*}
    \frac{v}{v(e,e,\ldots,e)} \circ m &=\frac{v}{v(e,e,\ldots,e)} \circ(u \circ m)\\
    &=\left(\frac{v u}{v(e,e,\ldots,e)}\right) \circ m=w \circ m=m.
\end{align*}

\item Now, the statement is proved for $v \in B^n(G)$ with $v(e,e,\ldots,e)=0$. Choose $u \in A^n(G)$ such that $u(e,e,\ldots,e)=1$. Then $((1-v)u)(e,e,\ldots,e)=1$ and hence by using (iii), $((1-v)u) \circ m=m$, and so $v \circ (u \circ m)=0$. And, again by using (iii),
$$
v \circ m= v \circ(u \circ m)=0=v(e,e,\ldots,e) m.
$$
Hence the result. \endproof
\end{enumerate}
\end{proof}
 Some properties of $M_{A^n(G)}$ and $M_{B^n(G)}$ are listed in the following lemma.
\begin{lem}\label{L1XCS5}
    Let $G$ be a locally compact group. Then we have the following.
     \begin{enumerate}[(i)]
         \item $M_{A^n(G)}$ and $M_{B^n(G)}$ are convex subsets of $A^n(G)$ and $B^n(G)$ respectively.
        \item $M_{A^n(G)}$ and $M_{B^n(G)}$ are abelian semigroups under pointwise multiplication.
        \item If $u \in M_{B^n(G)}$ and $m \in M$, then $u \circ m \in M$.
        \item If $\iota$ denotes the canonical inclusion of $A^n(G)$ into its second dual, then
        $$\iota\left(M_{A^n(G)}\right) \supset M.$$
\end{enumerate}
\end{lem}
\begin{proof}
    We will prove (iii) since the other parts are trivial. Let $v \in M_{B^n(G)}$ and $m \in M$. Then
$$
\begin{aligned}
1 & =\|m\|=\|v\|_{B^n(G)}\|m\| \geq\|v \circ m\| \\ & \geq \langle (I \otimes \ldots \otimes I), v \circ m\rangle\\ &=\langle v(e,e \ldots, e) (I \otimes \ldots \otimes I), m\rangle \langle (I \otimes \ldots \otimes I), m\rangle=\|m\|
\end{aligned}
$$
Thus $\|v \circ m\|=v \circ m(I \otimes \ldots \otimes I)=1$.
\end{proof}

\begin{prop}\label{L2XCS5}
There exists a mean $m \in M$ such that $u \circ m=m$ for all $u \in M_{A^n(G)}$.
\end{prop}
\begin{proof}
     Let $VN^n(G)^*$ be a Banach space and observe that $M$ is a non-empty weak$^*$ compact convex subset of $VN^n(G)^*$. Define the operators $T_u:M \rightarrow M$, $T_u(m) = u \circ m$ for all $u \in M_{A^n(G)}$, as weak* continuous affine operators (using Lemma \ref{L1XCS5}(iii)). By applying Markov- Kakutani Fixed Point Thoerem, there exists $m \in M$ such that $u \circ m= m$ for all $u \in M_{A^n(G)}.$
\end{proof}
\section{Invariant means on $VN^n(G)$}\label{sctn4CS5} 
The invariant means on $VN^n(G)$ is defined and their basic properties are studied in this section.

\begin{defn}
    \begin{itemize}
        \item A mean $m$ is called a topological invariant mean if 
        $$m(uT)= u(e,\ldots,e)m(T)\ \forall\ T\in VN^n(G),\ u \in A^n(G).$$
        \item $\operatorname{TIM}^n(\widehat{G})$ denotes the set of all topological invariant means on $VN^n(G)$.
        \item  $P^n_G$ denotes the set of elements in $A^n(G)$ with \\
        $u(e,e, \ldots, e)= \|u\|_{A^n(G)}=1$.
        \item A net $\{v_{\alpha}\} \subset P^n_G$ is called $w$ convergent to invariance if for all $v \in P^n_G$,
        $$ w-lim(v_{\alpha}v-v_{\alpha})=0,i.e., lim\langle v_{\alpha}v-v_{\alpha}, T \rangle =0$$
        for all $T \in VN^n(G).$
        \item A net $\{v_{\alpha}\} \subset P^n_G$ is called norm convergent to invariance if for all $v \in P^n_G$,
        $$ \|v_{\alpha}v-v_{\alpha}\|=0.$$
    \end{itemize}
\end{defn}

\begin{prop}\label{P2XCS5}
The set $\operatorname{TIM}^n(\widehat{G})$ is non-empty.
\end{prop}
\begin{proof}
    It is clear by a direct use of Theorem \ref{Thm1XCS5}.
\end{proof}

 Let us begin with some properties of an invariant functional on $VN^n(G)$.
    \begin{lem}
        $(VN^n(G))^*$ is a Banach algebra with the  operation $\circ$ where 
        $$\|m \circ p \| \leq \|m\| \|p\|\ \forall\ m,p \in (VN^n(G))^*. $$
    \end{lem}
    \begin{proof}
        The proof follows from the definition of the operation $\circ$.
    \end{proof}
\begin{lem}\label{L3XCS5}
\begin{enumerate}[(i)]
    \item  Let $m \in (VN^n(G))^* $, then $m$ is invariant iff $m= u \circ m\ \forall u \in P^n_G$.
    \item Let $m, p\in (VN^n(G))^* $ and if $m$ is invariant then so is $m \circ p.$
    \end{enumerate}
\end{lem}
\begin{proof}
\begin{enumerate}[(i)]
   \item Let $T \in VN^n(G)$ and $u \in P^n_G$, then 
    $$ \left\langle u \circ m, T \right\rangle = \left\langle u, m \circ T \right\rangle= \left\langle m \circ T, u \right\rangle = \left\langle m, uT\right\rangle.$$
    \item Observe that for any $u \in P^n_G$, $$ u \circ (m \circ p) = (u \circ m) \circ p = m \circ p.$$ 
    \end{enumerate}
\end{proof}
In the next theorem, we study some properties of invariant means.
\begin{thm}\label{Thm2XCS5}
\begin{enumerate}[(i)]
    \item If $G$ is discrete, then $\operatorname{TIM}^n(\widehat{G}) \cap A^n(G)$ is a singleton set.
    \item If $\operatorname{TIM}^n(\widehat{G}) \cap A^n(G) \neq \emptyset$ then $G$ is discrete.
    \item If $G$ is not discrete and $m \in \operatorname{TIM}^n(\widehat{G})$, then $m(T)=0$ for all $T \in \otimes^{n}_{h} C_{\rho}^{*}(G)$.
    \item If $G$ is non-discrete, $u \in M_{A^n(G)}$ and $m \in \operatorname{TIM}^n(\widehat{G})$ then $\|u-m\|=2$.
\end{enumerate}
\end{thm}
\begin{proof}
\begin{enumerate}[(i)]
    \item Observe that the function $\chi_e  \otimes\ldots \otimes \chi_e \in A^n(G)$ because $\chi_e \ast \tilde{\chi_e} = \chi_e $ and $G$ is discrete.
    Now, define a linear functional $m$ on $VN^n(G)$ as
    $$  m(T) = \left\langle T, \chi_e \otimes\ldots \otimes \chi_e  \right\rangle = \left \langle T(\chi_e \otimes\ldots \otimes \chi_e), \chi_e \otimes\ldots \otimes \chi_e \right \rangle. $$    From the above calculation, it is clear that $m$ is a mean.
    Let $T \in VN^n(G)$ and $u \in P^n_G$  then $u(\chi_e \otimes\ldots \otimes \chi_e)= \chi_e \otimes\ldots \otimes \chi_e.$
    Also,
    \begin{align*}
        m(uT) &= \left \langle uT, \chi_e \otimes\ldots \otimes \chi_e \right\rangle 
        = \left\langle T,u(\chi_e \otimes\ldots \otimes \chi_e)\right \rangle\\
        &=\left \langle T,\chi_e \otimes\ldots \otimes \chi_e\right \rangle = m(T).
    \end{align*}
     Hence, $m$ is an invariant mean on $VN^n(G)$. To prove the uniqueness of $m$, let us assume that $m^{'}$ is any invariant mean on $VN^n(G)$. Consider $T \in VN^n(G)$ and $v \in A^n(G)$, then
     \begin{align*}
        \left \langle (\chi_e \otimes\ldots \otimes \chi_e) T, v \right\rangle &= \left\langle T, (\chi_e \otimes\ldots \otimes \chi_e)v \right\rangle\\
         &= \left\langle T, v(e,\ldots, e)(\chi_e \otimes\ldots \otimes \chi_e)\right \rangle\\
         &= v(e,\ldots, e)\left\langle T, (\chi_e \otimes\ldots \otimes \chi_e)\right \rangle\\
         &= v(e,\ldots, e) m(T)\\
         &=\left \langle m(T)(I\otimes \ldots \otimes I),v \right\rangle,
     \end{align*}
     which shows $(\chi_e \otimes\ldots \otimes \chi_e) T=m(T)(I\otimes \ldots \otimes I).$ Since $m^{'}$ is invariant and $\chi_e \otimes\ldots \otimes \chi_e \in P^n_G,$
     \begin{align*}
        m^{'}(T) &= m^{'}(\chi_e \otimes\ldots \otimes \chi_e)T\\
        &= m^{'}(m(T)(I\otimes \ldots \otimes I))\\
        &=m(T).
     \end{align*}
     Hence the result.
\item Since $m \in A^n(G)$, then there exists a function $v \in A^n(G)$ such that 
    $$m(T) = \left\langle v,T \right\rangle\ \forall\ T \in VN^n(G) .$$
    Then $\left\langle v,T \right\rangle = m(T) = m(uT) = \left \langle uT, v \right\rangle = \left\langle T,uv \right\rangle\ \forall\ T \in VN^n(G)\ \forall\ u \in P^n_G. $
    It implies that $v=uv\ \forall\ u\in P^n_G.$
    Now, assume that the group $G$ is not discrete. Then we choose both $(x_1,\ldots, x_n) \in G^n$ and $u\in P^n_G$ such that $v(x_1,\ldots, x_n) \neq 0$ and  $u(x_1,\ldots, x_n) = 0$, which contradicts the relation $v=uv.$ Hence, the group $G$ is discrete.
    \item As $\otimes^n L^{1}(G)$ is dense in $\otimes_h^nC_{\rho}^{*}(G)$ in norm topology, so it is sufficient to prove that 
     $$m\left(L_{f_1} \otimes \ldots\otimes L_{f_n}\right)=0\ \forall\  f_1 \otimes \ldots \otimes f_n \in \otimes^n L^{1}(G),$$
     where $L_f$ denotes the bounded operator on $L_2(G)$ defined by $L_f x= f\ast x$, where $\ast$ is convolution. 
Let $f_1 \otimes \ldots \otimes f_n \in \otimes^n L^{1}(G), \epsilon>0$ and $V$ be a neighbourhood of $(e, \ldots, e)$ in $G^n$ such that
$$\int_{V}|f_1 \otimes \ldots \otimes f_n(g_1, \ldots, g_n)| dg_1 \ldots dg_n <\epsilon.$$ 
Now, choose $v \in P^n_G$ with supp $v \subset V$. Also, we may choose $\left\{v_{\alpha}\right\} \subset P^n_G$ such that w$^{*}$-$\lim v_{\alpha}=m$. Then $\left\{v_{\alpha}\right\}$ is $w$-convergent to invariance so that
\begin{align*}
&m\left(L_{f_1} \otimes \ldots\otimes L_{f_n}\right)= \lim \int f_1 \otimes \ldots \otimes f_n(g_1, \ldots, g_n) v_{\alpha}(g_1, \ldots, g_n) dg_1 \ldots dg_n\\
&=\lim \int f_1 \otimes \ldots \otimes f_n(g_1, \ldots, g_n)v(g_1, \ldots, g_n) v_{\alpha}(g_1, \ldots, g_n) dg_1 \ldots dg_n.
\end{align*}

But
\begin{align*}
&\left| \int f_1 \otimes \ldots \otimes f_n(g_1, \ldots, g_n)v(g_1, \ldots, g_n) v_{\alpha}(g_1, \ldots, g_n) dg_1 \ldots dg_n\right|\\
& \leq  \int \left|f_1 \otimes \ldots \otimes f_n(g_1, \ldots, g_n)v(g_1, \ldots, g_n) v_{\alpha}(g_1, \ldots, g_n)\right| dg_1 \ldots dg_n \\
& \leq \int_{V}|f_1 \otimes \ldots \otimes f_n(g_1, \ldots, g_n)|dg_1 \ldots dg_n <\epsilon \text { for all } \alpha
\end{align*}
Hence $m\left(L_{f_1} \otimes \ldots\otimes L_{f_n}\right)=0$.
\item Let $u \in M_{A^n(G)}$. By \cite[Section 4 and Theorem 5.1]{ToTu}, $u$ can be considered as a positive linear functional on $ \otimes_h^n C_{\rho}^{*}(G)$. By \cite[Corollary 3.5]{Taka}, $C_{\rho}^{*}(G)$ and hence $ \otimes^n_h C_{\rho}^{*}(G)$ has an approximate identity. Then for any $\epsilon>0$, there exists $S \in \otimes_h^n C_{\rho}^{*}G)$ such that $0 \leq S \leq I$ and $\langle S, u\rangle \geq 1-\epsilon$. Let $T=2 S-I \in VN^n(G)$. Then for $m \in \operatorname{TIM}^n(\widehat{G})$, since $\langle I, u-m\rangle=0$, we have
$$ \langle T, u-m\rangle=2\langle S, u-m\rangle \geq 2\langle S, u\rangle \geq 2(1-\epsilon)$$ Since $\epsilon>0$ is arbitrary and $\|T\| \leq 1$, we have $\|u-m\|=2$. \endproof
\end{enumerate}
\end{proof}
Let $H$ be a compact normal subgroup of $G$, then it is easy to see that there is a canonical inclusion of $A^n(G/H)$ inside $A^n(G)$, since  $A^n(G/H)$ can be viewed as a subalgebra of $A^n(G)$ consisting of all functions constant on cosets of $H^n$. Let $q_n$ denote the adjoint of above inclusion map from $VN^n(G)$ to $VN^n(G/H).$
In a similar way, $\otimes^{\sigma h}_nL^2(G/H)$ may be identified with a subspace $S$ of  $\otimes^{\sigma h}_nL^2(G)$.  Any $\tilde{m} \in VN^n(G/H)^*$ is said to be an invariant element if $\tilde{m}(\tilde{u} \widetilde{T})= \tilde{u}(\tilde{e},\tilde{e},...,\tilde{e})\tilde{m}(\widetilde{T})$ for all $\tilde{u} \in A^n(G/H)$ and $\widetilde{T} \in VN^n(G/H)$.
\begin{nota}
 For any $\widetilde{T} \in VN^n(G/H),$ the map
$\rho^n \widetilde{T}$ defines an operator on the tensor product $\otimes^{\sigma h}_nL^2(G)$ in the following manner:
$$ \rho^n \widetilde{T}(x_1,\ldots,x_n)= 
 \begin{cases} 
      \widetilde{T}(x_1,\ldots,x_n) & :(x_1,\ldots,x_n)\in S \\
      0 & :(x_1,\ldots,x_n) \notin S. 
   \end{cases}
$$

For $v \in A^n(G)$, $f(\widetilde{T}) = \langle v, \rho^n \widetilde{T} \rangle$ is an ultraweakly continuous linear functional on $VN^n(G/H)$ so that $f(\widetilde{T}) = \langle \tilde{v}, \widetilde{T} \rangle$ for some uniquely defined $\tilde{v} \in A^n(G/H).$
\end{nota}
\begin{prop}\label{P3XCS5}
    Let $H$ be a normal compact subgroup of $G$ and $m \in \operatorname{TIM}^n(\widehat{G})$ and $\tilde{p}$ be an invariant element in $(VN^n(G/H))^*.$ Consider $p= q_n^* \tilde{p}$, then $m \circ p$ is invariant  and $\|m \circ p\|=  \|\tilde{p}\|.$
\end{prop}
\begin{proof}
    Using Lemma \ref{L3XCS5} (ii), it is clear that $m \circ p$ is invariant. Also, 
    $$\|m \circ p \| \leq \|m\| \|p\| = \|q^*_n(\tilde{p})\| \leq \|\tilde{p}\|. $$
     Now, fix $\epsilon >0,$ and choose $\widetilde{T} \in VN^n(G/H)$ with $\|\widetilde{T}\| =1$ and 
     $$ |\left\langle \tilde{p}, \widetilde{T} \right\rangle | \geq \| \tilde{p}\| - \epsilon.$$
     Let $v \in A^n(G)$ and $\tilde{u} \in A^n(G/H)$ which is identified as a element of $A^n(G)$, then
     \begin{align*}
         \left\langle q_n(v \rho^n \widetilde{T}), \tilde{u} \right\rangle &= \left\langle v \rho^n \widetilde{T}, \tilde{u} \right\rangle = \left\langle \rho^n \widetilde{T}, v\tilde{u} \right\rangle \\
         &= \left\langle\tilde{u} \rho^n \widetilde{T}, v \right\rangle = \left\langle\tilde{u} \widetilde{T}, \tilde{v} \right\rangle\\
         &= \left\langle\tilde{v} \widetilde{T}, \tilde{u} \right\rangle
     \end{align*}
     It implies that $q_n(v \rho^n \widetilde{T})=\tilde{v} \widetilde{T}.$
     Now, 
     \begin{align*}
        \left \langle p \circ \rho^n \widetilde{T}, v\right\rangle &= \left\langle p,v \rho^n \widetilde{T} \right\rangle = \left\langle q^*_n \tilde{p}, v \rho^n \widetilde{T} \right\rangle\\
         &= \left\langle \tilde{p}, q_n^*(v \rho^n \widetilde{T})\right \rangle= \left\langle \tilde{p}, \tilde{v} \widetilde{T} \right\rangle\\
         &= \tilde{v}(\tilde{e},\ldots, \tilde{e} )\left\langle \tilde{p}, \widetilde{T} \right\rangle .
     \end{align*}
     Hence, $p \circ \rho^n \widetilde{T}= \left\langle \tilde{p}, \widetilde{T} \right\rangle(I \otimes \ldots \otimes I)$ and we have
     \begin{align*}
         |\left\langle m \circ p , \rho^n \widetilde{T} \right\rangle| &= |\left\langle m , p \circ \rho^n \widetilde{T} \right\rangle|\\
         &|\left\langle m , \left\langle \tilde{p}, \widetilde{T} \right\rangle(I \otimes \ldots \otimes I) \right\rangle|\\
         &= |\left\langle \tilde{p}, \widetilde{T} \right\rangle| \geq \|\tilde{p}\| - \epsilon.
     \end{align*}
       It shows that $\|m \circ p \| \geq \|\tilde{p}\|$, and this completes the proof of the result.
\end{proof}

\begin{cor}\label{C1XCS5}
    If $H$ is compact normal and $VN^n(G/H)^*$ has more than one invariant mean, then so does $(VN^n(G))^*.$
\end{cor}
\begin{proof}
    Let $m$ be the invariant mean on $VN^n(G).$
    Define a map 
    $$\sigma_n : (VN^n(G/H))^* \rightarrow (VN^n(G))^* \text{ as }
    \sigma_n (\tilde{p}) = m \circ q_n^* \tilde{p} $$

    Observe that $\|\sigma_n (\tilde{p})\|= \|m \circ q_n^* \tilde{p}\|= \|m \circ p\| =\|\tilde{p}\|$
    (From Proposition \ref{P3XCS5}). Hence, the above map is an isometry  from the subspace of invariant elements in $(VN^n(G/H))^*$ into the subspace of invariant elements in $(VN^n(G))^*.$
     As a result, if $\tilde{p_1},\tilde{p_2}$ are two distinct invariant means in $(VN^n(G/H))^*,$ then $\sigma_n(\tilde{p_1}),\sigma_n(\tilde{p_2})$ are distinct invariant means in $(VN^n(G))^*.$
\end{proof}

\begin{prop}\label{P4XCS5}
    Let $S$ be an open subgroup of $G$. If there exist two distinct invariant means on $VN^n(S),$ then there are two distinct invariant means on $VN^n(G).$ 
\end{prop}

\begin{proof}
  Suppose $r_n$ denotes the restriction homomorphism from $A^n(G)$ to $A^n(S)$. From  \cite[Lemma 3.1]{Kanu}, the map $r_n$ is bounded, and because $S$ is open, it is also onto. We claim that if $m^{'}$ is an invariant mean in $VN^n(S)$, then there exists an invariant mean $m$ in $VN^n(G),$ such that $m^{'}= m|_{(r_n)^*(VN^n(S)}$.
  If $m^{'}$ is an invariant mean on $VN^n(S)$, then there exists a net $\{w_\alpha^{'}\}$ of normalized functions in $A^n(S)$ with its w$^*$- limit equal to $m^{'}$. Let $U$ be any neighborhood of $(e,e,\ldots, e)$ in $G^n$ and we choose $v_{U} \in P^n_G$ with that supp $v_{U} \subset U$. Now consider $w_{\alpha, U}^{'}=w_{\alpha}^{'} v_{U}|_S$ and observe that for any $T \in VN^n(S)$ and for all $U$,

$$
\lim \left\langle w_{\alpha, U}^{'}, T\right\rangle=\lim \langle w_{\alpha}^{'}, v_{U}|_S \cdot T\rangle= \left\langle m^{'}, v_{U}|_S \cdot T\right\rangle=\left\langle m^{\prime}, T\right\rangle
$$

Choose $w_{\alpha} \in P^n_G$ with the property $r_n w_{\alpha}=w_{\alpha}^{'}$ and now define $w_{\alpha, U}=w_{\alpha} v_{U}$. If $m$ represents the $w^{*}$-limit point of $\left\{w_{\alpha, U}\right\}$, then $m$ is an invariant mean on $VN^n(G)$, since $\{w_{\alpha, U}^{'}\}$ is w$^*$- convergent to invariance. Next, we need prove that $$m|_{(r_n)^{*} (VN^n(S))}=m^{'}.$$

For $T \in VN^n(S), \epsilon>0$, we may choose $\alpha_{0}$, such that for $\alpha \geq \alpha_{0}$ and all $U$ $$|\left\langle w_{\alpha, U}^{'}, T\right\rangle -\left\langle m^{'}, T\right\rangle|<\epsilon / 2.$$ 
Now, choose $(\alpha, U) \geq\left(\alpha_{0}, U\right)$ such that $\left|\left\langle m, (r_n)^{*} T\right\rangle-\left\langle w_{\alpha, U}, (r_n)^{*} T\right\rangle\right|<$ $\epsilon / 2$. Then
\begin{align*}
& \left|\left\langle m, (r_n)^{*} T\right\rangle-\left\langle m^{'}, T\right\rangle\right| \\
&\leq\left|\left\langle m, (r_n)^{*} T\right\rangle-\left\langle w_{\alpha, U}, (r_n)^{*} T\right\rangle\right|+\left|\left\langle w_{\alpha, U}, (r_n)^{*} T\right\rangle-\left\langle m^{'}, T\right\rangle\right| \\
& <\epsilon / 2+\left|\left\langle r_n w_{\alpha, U}^{'}, T\right\rangle-\left\langle m^{\prime}, T\right\rangle\right|\\
&=\epsilon / 2+\left|\left\langle w_{\alpha, U}^{'}, T\right\rangle-\left\langle m^{'}, T\right\rangle\right|<\epsilon
\end{align*}
Hence $\left\langle m, (\rho^n)^{*} T\right\rangle=\left\langle m^{'}, T\right\rangle$.
\end{proof}
\begin{thm}
    Let $G$ be a second countable group. If $VN^n(G)$ admits a unique topological invariant mean, then $G$ is discrete.
\end{thm}
\begin{proof}

    Let us denote $\mathscr{U}$ by a neighbourhood basis of $(e,e,\ldots, e)$ with the property that each element of $\mathscr{U}$ is a compact set. Without loss of generality, we may assume that $\mathscr{U}$ is countable, since $G$ is second countable. Hence assume that $\mathscr{U}$ is a sequence $\{\widetilde{U}_{n}\}$ such that $\widetilde{U}_{n} \rightarrow\{(e,e,\ldots,e)\}$. For each $n \in \mathbb{N}$, fix $v_{n} \in M_{A^n(G)}$ with $\operatorname{supp}\left(v_{n}\right) \subseteq \widetilde{U}_{n}$.

Let $v \in M_{A^n(G)}$ and $\epsilon>0$.The set of compactly supported elements in $M_{A^n(G)}$ is dense in $M_{A^n(G)}$ since compactly supported elements of $M_{A(G)}$ are dense in $M_{A(G)}$(\cite[Corollory 2.3.5]{KL2}). Hence there exists $v^{'} \in M_{A^n(G)}$ with compact support such that $\left\|v-v^{'}\right\|<\epsilon / 2$. By regularity of $A^n(G)$( follows on the same steps as \cite[Proposition 3.2]{alaghmandan}), there exists $u \in A^n(G)$ such that $u$ is 1 on $\operatorname{supp}\left(v^{'}\right)$. Since $v^{'}(e,e,\ldots,e)=1$, we have $(e,e,\ldots,e) \in \operatorname{supp}\left(v^{'}\right)$, and $\left(v^{'}-u\right)(e,e,\ldots,e)=0$.

From \cite{Kanu}, $\{(e,e,\ldots,e)\}$ is a set of spectral synthesis and hence there exists $w \in A^n(G)$ such that $\left\|v^{'}-u-w\right\|<\epsilon / 2$ and $w(\widetilde{W})=0$ for some neighborhood $\widetilde{W}$ of $(e,e,\ldots,e)$. Now, for any $n \in \mathbb{N}$ such that $\widetilde{U}_{n} \subset \widetilde{W} \cap \operatorname{supp}\left(v^{'}\right)$, we can observe that $v_n u=v_{n}$ and $v_{n} w=0$. And, by using $\epsilon / 2$ argument, we have $$\left\|v v_{n}-v_{n}\right\|= \|vv_n-v^{'}v_n+v^{'}v_n-uv_n+uv_n-wv_n+wv_n-v_n\|$$ $$\leq \|v-v^{'}\|\|v_n\|+\|v^{'}-u-w\|\|v_n\|+\|uv_n+wv_n-v_n\|<\epsilon/2 +\epsilon/2+0=\epsilon$$.

It follows that every weak* accumulation point of $\left\{v_{n}\right\}$ in $A^n(G)^{* *}$ is a topological invariant mean. By the assumption of uniqueness of topological invariant mean and as the set of topological invariant means on $VN^n(G)$ is nonempty, let $m$ denote the unique topological invariant mean on $VN^n(G)$. Since $A^n(G)$ is the predual of the von Neumann algebra $VN^n(G)$, it is weakly sequentially complete. Thus $\left\{v_{n}\right\}$ converges to $m$ weakly in $A^n(G)$, which means that $m \in A^n(G)$. Hence by Theorem \ref{Thm2XCS5} (ii), it follows that $G$ is discrete.
\end{proof}

\section{Open Subgroups and Invariant means} \label{sctn5CS5}
We initially provide a list of functorial properties of $A^n(G)$ in accordance with the principles of \cite{KL2}. We omit the proofs, as they are similar to those in \cite{KL2}.
\begin{lem}\label{L3CS5}
Let $H$ be a closed subgroup of $G^n$. For $u \in A^n(H)$ let $u^{\circ}$ denote the extension function on $G$ that is $u$ on $H^n$ and vanishes outside $H^n$.\\
(i) If $H$ is open, then $ u \mapsto u^{\circ}$ is an isometric isomorphism of $A^n(H)$ onto $A^n(G)^{\circ}=\left\{u^{\circ}: u \in A^n(H)\right\}$\\
(ii) The restriction map $ A^n(G) \mapsto A^n(H)$ is a surjective contractive homomorphism.
\end{lem}
\begin{nota}
    Let $H$ be an open subgroup of the locally compact group $G$. We denote $r_n$ as the restriction map from $A^n(G)$ to $A^n(H)$ and $e_n$ as the extension map from $A^n(H)$ to $A^n(G).$ Observe that $r_ne_n$ is actually the identity of $A^n(H)$.
\end{nota} 
\begin{lem}\label{L5CS5}
    For $u \in A^n(G)$ and $T \in VN^n(G)$, we have
$$
e_{n}^{*}(u \cdot T)=r_{n}(u) \cdot e_{n}^{*}(T)
$$
\end{lem}
\begin{proof}
    For any $v \in A^n(H)$,

$$
\begin{aligned}
\left\langle v, e_{n}^{*}(u \cdot T)\right\rangle & =\left\langle e_{n}(v), u \cdot T\right\rangle=\left\langle u e_{n}(v), T\right\rangle \\
& =\left\langle e_{n}\left(r_{n}(u) v\right), T\right\rangle=\left\langle r_{n}(u) v, e_{n}^{*}(T)\right\rangle \\
& =\left\langle v, r_{n}(u) \cdot e_{n}^{*}(T)\right\rangle
\end{aligned}
$$
\end{proof}
\begin{lem}\label{L7CS5}
    For $u \in A^n(G)$ and $T \in VN^n(G)$, we have
$$
r_{n}^{*}(u \cdot T)=e_{n}(u) \cdot r_{n}^{*}(T)
$$
\end{lem}
\begin{proof}
    For any $v \in A^n(G)$,

$$
\begin{aligned}
\left\langle v, r_{n}^{*}(u \cdot T)\right\rangle & =\left\langle r_{n}(v), u \cdot T\right\rangle=\left\langle u r_{n}(v), T\right\rangle \\
& =\left\langle r_{n}\left(e_{n}(u) v\right), T\right\rangle=\left\langle e_{n}(u) v, r_{n}^{*}(T)\right\rangle \\
& =\left\langle v, e_{n}(u) \cdot r_{n}^{*}(T)\right\rangle 
\end{aligned}
$$
\end{proof}
\begin{lem}\label{L6CS5 sec5}
    For $v \in A^n(H)$ and $T \in VN^n(G)$, we have
$$
r_{n}^{*}\left(v \cdot e_{n}^{*}(T)\right)=e_{n}(v) \cdot T
$$
\end{lem}

\begin{proof}
    For any $u \in A^n(G)$,

$$
\begin{aligned}
\left\langle u, r_{n}^{*}\left(v \cdot e_{n}^{*}(T)\right)\right\rangle & =\left\langle r_{n}(u), v \cdot e_{n}^{*}(T)\right\rangle=\left\langle v r_{n}(u), e_{n}^{*}(T)\right\rangle \\
& \left.=\left\langle e_{n}\left(v \cdot r_{n}(u)\right), T\right)\right\rangle=\left\langle e_{n}(v) u, T\right\rangle \\
& =\left\langle u, e_{n}(v) \cdot T\right\rangle
\end{aligned}
$$
\end{proof}
\begin{lem}
   The second adjoint of the extension map  $e_{n}^{* *}: A^n(H)^{* *} \rightarrow A^n(G)^{* *}$ is an isometry.
\end{lem}
\begin{proof}
     The restriction map $r_{n}: A^n(G) \rightarrow A^n(H)$ is a contraction by Lemma \ref{L3CS5} and hence $\left\|r_{n}^{* *}\right\|=\left\|r_{n}\right\| \leq 1$. Also, $r_{n}^{* *} e_{n}^{* *}$ is the identity map on $A^n(H)^{* *}$ and we have $\|\phi\|=\left\|r_{n}^{* *} e_{n}^{* *}(\phi)\right\|$ for any $\phi \in A^n(G)^{* *}$. Now, let $\left\|e_{n}^{* *}(\phi)\right\|<\|\phi\|$ for some $\phi \in A^n(G)^{* *}$. Then
$$
\|\phi\|=\left\|r_{n}^{* *} e_{n}^{* *}(\phi)\right\| \leq\left\|e_{n}^{* *}(\phi)\right\|<\|\phi\|
$$
which is actually a contradiction.
\end{proof}
We will now proceed to show the key results of this section.
\begin{thm}\label{Thm2CS5}
    Let $H$ be an open subgroup of $G$, \\
    (i) If $m \in \operatorname{TIM}^n(\widehat{H})$ then $e_n^{* *}(m) \in \operatorname{TIM}^n(\widehat{G}).$\\
    (ii) If $m \in \operatorname{TIM}^n(\widehat{G})$, then $\exists\ \widetilde{m} \in  \operatorname{TIM}^n(\widehat{H})$ such that  $e_n^{* *}(\widetilde{m})= m$
\end{thm}
\begin{proof}
   (i) Suppose that $m \in \operatorname{TIM}^n(\widehat{H})$ and let $u \in A^n(G)$ and $T \in VN^n(G).$ Then by Lemma \ref{L5CS5},
    $$
\begin{aligned}
\left\langle u \cdot T, e_{n}^{* *}(m)\right\rangle & =\left\langle e_{n}^{*}(u \cdot T), m\right\rangle=\left\langle r_{n}(u) \cdot e_{n}^{*}(T), m\right\rangle \\
& =r_{n}(u)(e,e,\ldots,e)\left\langle e_{n}^{*}(T, m)\right\rangle\\
&=u(e,e,\ldots,e)\left\langle T, e_{n}^{* *}(m)\right\rangle.
\end{aligned}
$$
It implies that $e_n^{* *}(m) \in \operatorname{TIM}^n(\widehat{G})$ and hence the result.\\
(ii) Suppose that $m \in \operatorname{TIM}^n(\widehat{G})$, then we claim that $$\widetilde{m} = r_n^{* *}(m).$$ 
We first prove that $e_n^{* *}(r_n^{* *}(m))=m.$ Let $u \in A^n(H)$ be such that $u(e,e,\ldots,e)=1$ and $T \in VN^n(G)$. Using Lemma \ref{L6CS5 sec5}, we get 
$$
\begin{aligned}
\langle T, m\rangle & =u(e,e,\ldots,e)\langle T, m\rangle=e_{n}(u)(e,e,\ldots,e)\langle T, m\rangle \\
& =\left\langle e_{n}(u) \cdot T, m\right\rangle=\left\langle r_{n}^{*}\left(u \cdot e_{n}^{*}(T)\right), m\right\rangle \\
& =\left\langle u \cdot e_{n}^{*}(T), r_{n}^{* *}(m)\right\rangle=\left\langle e_{n}^{*}(T), r_{n}^{* *}(m)\right\rangle \\
& =\left\langle T, e_{n}^{* *} r_{n}^{* *}(m)\right\rangle
\end{aligned}
$$
Next, we will prove that $\widetilde{m}$ is an invariant mean on $VN^n(H)$. Let $u \in A^n(H)$ and $T \in VN^n(H)$, then by using Lemma \ref{L7CS5}, we have
$$
\begin{aligned}
\left\langle u \cdot T, r_{n}^{* *}(m)\right\rangle & =\left\langle r_{n}^{*}(u \cdot T), m\right\rangle=\left\langle e_{n}(u) \cdot r_{n}^{*}(T), m\right\rangle \\
& =e_{n}(u)(e,e,\ldots,e)\left\langle r_{n}^{*}(T), m\right\rangle=u(e,e,\ldots,e)\left\langle T, r_{n}^{* *}(m)\right\rangle.
\end{aligned}
$$
Hence the result.
\end{proof}
\begin{nota}
    Let $\#S$ denote the cardinality of a set $S$.
\end{nota}
    \begin{cor}\label{C2CS5}
        Let $H$ be an open subgroup of the locally compact group $G$, then\\ 
        (i) $\# \operatorname{TIM}^n(\widehat{H}) = \#\operatorname{TIM}^n(\widehat{G})$. \\
        (ii) $\operatorname{TIM}^n(\widehat{H})$ is separable iff $\operatorname{TIM}^n(\widehat{G})$ is separable.
    \end{cor}
    \begin{proof}
       Since the map $e_n^{**}$ is injective, the proof follows directly from Theorem \ref{Thm2CS5}.
    \end{proof}
    
\begin{cor}
    If $G$ is discrete, then there exists a unique invariant mean on $VN^n(G)$.
\end{cor}
\begin{proof}
    Let $H = \{e\}$, then by Corollary \ref{C2CS5}(i), we get the desired result.
\end{proof}
\begin{thm}\label{Thm3CS5}
    Let $H$ be a closed subgroup of $G$, then 
    $$ r_n^{* *}(\operatorname{TIM}^n(\widehat{G})) = \operatorname{TIM}^n(\widehat{H}).$$
\end{thm}

\begin{proof}
 Let $m \in r_{n}^{* *}(\operatorname{TIM}^n(\widehat{G}))$, then  $\exists\ \widetilde{m} \in \operatorname{TIM}^n(\widehat{G})$ such that $r_{n}^{* *}(\widetilde{m})=m$. Now, let $u \in A^n(H)$ and $T \in VN^n(H)$. Choose $v \in A^n(G)$ such that $r_{n}(v)=u$. Then 
 $$\langle u \cdot T, m\rangle= u(e,e,\ldots,e)\langle T, m\rangle$$ 
 follows exactly as in Theorem \ref{Thm2CS5}. Thus $r_{n}^{* *}(\operatorname{TIM}^n(\widehat{G})) \subseteq \operatorname{TIM}^n(\widehat{H}).$

For the other inclusion, consider an element $m$ belonging to the set $\operatorname{TIM}^n(\widehat{H})$. We claim the existence of a mean $\widetilde{m}$ on $VN^n(G)$ such that $r_{n}(\widetilde{m})=m$. Given that $m \in \operatorname{TIM}^n(\widehat{H})$, we may conclude that there exists a net $\left\{u_{\alpha}\right\}$ in $M_{A^n(H)}$ that converges to $m$ in the weak ${ }^{*}$-topology. Since $r_{n}\left(M_{A^n(G)}\right)=M_{A^n(H)}$, we might conclude that there is a net $\widetilde{u}_{\alpha} \subset M_{A^n(G)}$ such that $r_{n}\left(\widetilde{u}_{\alpha}\right)=u_{\alpha}$. According to the Banach-Alaoglu theorem, the set $\left\{\widetilde{u}_{\alpha}\right\}$ contains a subnet $\left\{\widetilde{u}_{\alpha_{\beta}}\right\}$ that converges weak${ }^{*}$ to $\widetilde{m}$. Since $r_{n}^{* *}$ is weak${ }^{*}$-weak${ }^{*}$-continuous, it follows that $r_{n}^{* *}(\widetilde{m})=m$.

Now, we claim that for any $u \in M_{A^n(G)}, r_{n}^{* *}(u \cdot \widetilde{m})=m$. For any $T \in VN^n(H)$, we have
$$
\begin{aligned}
\left\langle T, r_{n}^{* *}(u \cdot \widetilde{m})\right\rangle & =\left\langle r_{n}^{*}(T), u \cdot \widetilde{m}\right\rangle=\left\langle u \cdot r_{n}^{*}(T), \widetilde{m}\right\rangle \\
& =\left\langle r_{n}^{*}\left(r_{n}(u) \cdot T\right), \widetilde{m}\right\rangle=\left\langle r_{n}(u) \cdot T, r_{n}^{* *}(\widetilde{m})\right\rangle \\
& =\left\langle r_{n}(u) \cdot T, m\right\rangle=\langle T, m\rangle
\end{aligned}
$$
Hence the claim.

Now, consider the collection
$$
\mathcal{C}=\left\{\widetilde{m}: \widetilde{m} \text { is a mean on } VN^n(G) \text { and } r_{n}^{* *}(\widetilde{m})=m\right\}
$$
We have already seen that the collection $\mathcal{C}$ is non-empty. Also, it is clear that $\mathcal{C}$ is weak -${*}$ compact. Observe that for any $u \in M_{A^n(G)}$, the mapping $\widetilde{m} \mapsto u \cdot \widetilde{m}$ maps $\mathcal{C}$ to itself . Hence, the remaining proof will follow from the Markov-Kakutani fixed-point theorem.
\end{proof}
\section{Invariant means on the dual of $A^n_0(G)$}\label{sctn6CS5}
In the following section, we introduce a new algebra, denoted by $A^n_{0}(G)$, and examine the invariant means on the dual of this algebra. Let $M_nA(G)$ denote the collection of all multipliers of $A^n(G)$, ; that is, the set of all $\phi \in L^{\infty}(G^n)$ such that, for every $f \in A(G)$ we have $\phi\theta f \in A^n(G)$, where $\theta f(x_1,x_2,...,x_n)=f(x_1x_2...x_n).$
If the map $f \rightarrow \phi\theta f$ is completely bounded then we say that $\phi$ is a cb-multiplier of $A^n(G)$, and the collection of all cb-multiplier of $A^n(G)$ is denoted by $M_n^{cb} A(G)$. The cb-multiplier norm of $\phi$ is defined by
$$
\|\phi\|_{M_n^{c b} A(G)}=\left\|M_\phi\right\|_{c b},
$$
where
$$
M_\phi: A(G) \rightarrow A^n(G), \quad M_\phi(f)=\phi\theta f,
$$
and $\left\|M_\phi\right\|_{c b}$ denotes the completely bounded norm of the map $M_\phi$. Moreover, $B^n(G)\subset M_n^{cb} A(G)$. For more details about the cb-multiplier norm, one may refer to \cite{ToTu}.
\begin{nota}
   \begin{itemize}
       \item The closure of $A^n(G)$ in the cb-multiplier norm will be denoted as $A^n_{0}(G)$. 
       \item $\iota$ denotes the inclusion map from $A^n(G)$ to $A^n_{0}(G)$.
   \end{itemize}
\end{nota}
\begin{lem}\label{L6CS5}
   (i) For $u \in A^n(G)$ and $T \in VN^n(G), u \cdot T \in \iota^{*}\left(A_{0}^n(G)^{*}\right)$.\\
(ii) For $u \in A_{0}^n(G)$ and $T \in A_{0}^n(G)^{*}, \iota^{*}(u \cdot T)= u \cdot \iota^{*}(T)$. 
\end{lem}
\begin{proof} (i) Let us define the map 
$f_{u, T}: A_{0}^n(G) \rightarrow \mathbb{C}$ as $f_{u, T}(v)=\langle u v, T\rangle$. Observe that $f_{u, T}$ defines a bounded linear functional on $A_{0}^n(G)$ and $\left\|f_{u, T}\right\| \leq\|u\|_{A(^nG)}\|T\|_{VN^n(G)}$. Now, it is clear that $f_{u, T}$, when restricted to $A^n(G)$, equals $u \cdot T$, so it follows that $u \cdot T=\iota^{*}\left(f_{u, T}\right)$.\\
(ii) For any $v \in A^n(G)$, we have
$$
\begin{aligned}
\left\langle v, \iota^{*}(u \cdot T)\right\rangle & =\langle\iota(v), u \cdot T\rangle=\langle u \iota(v), T\rangle \\
& =\langle\iota(u v), T\rangle=\left\langle u v, \iota^{*}(T)\right\rangle \\
& =\left\langle v, u \cdot \iota^{*}(T)\right\rangle
\end{aligned}
$$
\end{proof}
\begin{defn}
\begin{enumerate}[(i)]
    \item  A linear functional $m$ on $A_{0}^n(G)^{*}$ is called a mean if $\|m\|=$ $m(I)=1$, where $I$ denotes the identity in $A_{0}^n(G)^{*}$.
   \item A mean $m$ on $A_{0}^n(G)^{*}$ is said to be  invariant if $u \circ m= u(e,e,\ldots,e) m$ for all $u \in A_{0}^n(G)$ i.e., 
   $$\langle T, u \circ m\rangle=\langle u . T, m\rangle= u(e,e,\ldots,e)\langle T, m\rangle\ \forall\ T \in A_{0}^n(G)^{*}\ \forall\ u \in A_{0}^n(G).$$
\end{enumerate}
\end{defn}
\begin{prop}
  If $(\operatorname{TIM}^n(\widehat{G}))_0$ denotes the set of all invariant means on $A_{0}^n(G),$ then $\iota^{* *}(\operatorname{TIM}^n(\widehat{G})) \subseteq (\operatorname{TIM}^n(\widehat{G}))_0$.
\end{prop}
\begin{proof}
    Let $m \in \operatorname{TIM}^n(\widehat{G})$. Let $u \in A_{0}^n(G)$ and $T \in A_{0}^n(G)^{*}$. As $u \in A_{0}^n(G)$, then there exists a sequence $\left\{u_{n}\right\} \subset A^n(G)$ such that $u_{n} \rightarrow u$ in the cb-multiplier norm. Hence $u_{n}(e,e,\ldots,e) \rightarrow u(e,e,\ldots,e)$. Now, for any $v \in A_{0}^n(G)$
$$
\begin{aligned}
\left|\left\langle v, u_{n} \cdot T- u \cdot T\right\rangle\right| & =\left|\left\langle v,\left(u_{n}-u\right) \cdot T\right\rangle\right|=\left|\left\langle\left(u_{n}-u\right) v, T\right\rangle\right| \\
& \leq\left\|u_{n}-u\right\|_{A_{0}^n(G)}\|v\|_{A_{0}^n(G)}\|T\|_{A_{0}^n(G)^{*}} .
\end{aligned}
$$
So, it implies that $u_{n} \cdot T \rightarrow u \cdot T$ in $A_{0}^n(G)^{*}$. And, we have
$$
\begin{aligned}
2\left\langle u \cdot T, \iota^{* *}(m)\right\rangle & =\left\langle\lim _{n \rightarrow \infty} u_{n} \cdot T, \iota^{* *}(m)\right\rangle=\lim_{n \rightarrow \infty}\left\langle u_{n} \cdot T, \iota^{* *}(m)\right\rangle \\
& =\lim_{n \rightarrow \infty}\left\langle\iota^{*}\left(u_{n} \cdot T\right), m\right\rangle \\
& =\lim_{n \rightarrow \infty}\left\langle u_{n} \cdot \iota^{*}(T), m\right\rangle \quad \quad \text { (by Lemma \ref{L6CS5}) } \\
& =\lim_{n \rightarrow \infty} u_{n}(e,e,\ldots,e)\left\langle\iota^{*}(T), m\right\rangle\\
&=u(e,e,\ldots,e)\left\langle\iota^{*}(T), m\right\rangle \\
& =u(e,e,\ldots,e)\left\langle T, \iota^{* *}(m)\right\rangle.
\end{aligned}
$$
\end{proof}
\begin{thm}\label{Thm4CS5}
    Let $G$ be a locally compact group. Then
$$
\iota^{* *}: \operatorname{TIM}^n(\widehat{G}) \rightarrow (\operatorname{TIM}^n(\widehat{G}))_0 
$$
is a bijection.
\end{thm}
\begin{proof}
First, we shall show the injectivity of the map $\iota^{* *}: \operatorname{TIM}^n(\widehat{G}) \rightarrow (\operatorname{TIM}^n(\widehat{G}))_0$. Let $m_{1}, m_{2} \in$ $\operatorname{TIM}^n(\widehat{G})$ with $m_{1} \neq m_{2}$. Then there exists $T \in VN^n(G)$ such that $m_{1}(T) \neq m_{2}(T)$. Now, choose an element $u_{0} \in A^n(G)$ with $u_{0}(e,e,\ldots,e)=1$. Then, we have
$$
\left\langle m_{1}, u_{0} \cdot T\right\rangle=\left\langle m_{1}, T\right\rangle \neq\left\langle m_{2}, T\right\rangle=\left\langle m_{2}, u_{0} \cdot T\right\rangle
$$
Using  Lemma \ref{L6CS5}, $u \cdot T \in A_{0}^n(G)^{*}$. And, we have
$$
\begin{array}{rlrl}
\left\langle u_{0} \cdot T, \iota^{* *}\left(m_{1}\right)\right\rangle & =\left\langle\iota^{*}\left(u_{0} \cdot T\right), m_{1}\right\rangle=\left\langle u_{0} \cdot T, m_{1}\right\rangle  \\
& \neq\left\langle u_{0} \cdot T, m_{2}\right\rangle=\left\langle\iota^{*}\left(u_{0} \cdot T\right), m_{2}\right\rangle  \\
& =\left\langle u_{0} \cdot T, \iota^{* *}\left(m_{2}\right)\right\rangle
\end{array}
$$
Now, for surjectivity, let $m \in (\operatorname{TIM}^n(\widehat{G}))_0$ and define $\widetilde{m} \in VN^n(G)^{*}$ as $\langle T, \widetilde{m}\rangle:=\langle u . T, m\rangle$, for all $T \in VN^n(G)$, where $u \in A^n(G)$ with $u_{0}(e,e,\ldots,e)=1$. It can be observed that $\widetilde{m}$ is actually independent of the choice of $u_{0}$ and hence it implies that $\widetilde{m} \in (\operatorname{TIM}^n(\widehat{G}))_0$. Now for any $T \in A_{0}^n(G)^{*}$, we have
$$
\left\langle T, \iota^{* *}(\widetilde{m})\right\rangle=\left\langle\iota^{*}(T), \widetilde{m}\right\rangle=\left\langle u_{0} \cdot \iota^{*}(T), m\right\rangle=\left\langle\iota^{*}(T), m\right\rangle=\langle T, m\rangle .
$$
Hence $\iota^{* *}(\widetilde{m})=m$.
\end{proof}
\begin{lem}
    Let $H$ be a closed subgroup of $G$. Then the restriction map $r_{n}$ is a contraction from $A_{0}^n(G)$ into $A_{0}^n(H)$.
\end{lem}
\begin{cor}
    Let $H$ be a closed subgroup of $G$. Then
    $$r_{n}^{* *}\left((\operatorname{TIM}^n(\widehat{G}))_0\right)=(\operatorname{TIM}^n(\widehat{H}))_0.$$
\end{cor}
\begin{proof}
The proof follows directly by using Theorem \ref{Thm3CS5} and Theorem \ref{Thm4CS5}.
\end{proof}

\end{document}